\documentclass[11pt,leqn0]{amsart}
\usepackage{epsf,amsmath,amsfonts,amsthm,amssymb,epsfig,enumerate,setspace}

\newtheorem*{thm5.1}{Theorem 5.1}
\newtheorem*{thm6.1}{Theorem 6.1}
\newtheorem*{lem6.2}{Lemma 6.2}
\newtheorem*{thm1}{Theorem 1}
\newtheorem*{thm2}{Theorem 2}
\newtheorem*{lem}{Lemma}

\newcommand{\E}{\mathbb{E}}

\newcommand{\prob}{\mathbb{P}}

\newcommand{\R}{\mathbb{R}}

\newcommand{\N}{\mathbb{N}}

\newcommand{\distrib}{\overset{\mbox{$d$}}{=}}

\title{A dynamical characterization of Poisson-Dirichlet distributions}
\author{Louis-Pierre Arguin}\thanks{The author is thankful to Michael Aizenman for introducing him to the subject and for plenty of insightful discussions. The author thanks also the referee for valuable remarks on the earlier version of the article. This work was supported in part by NSERC and FQRNT postgraduate fellowships and NSF grant DMS 0602360}
\address{Department of Mathematics, Princeton University, Princeton NJ 08544}
\email{larguin@math.princeton.edu}
\date{August 17th, 2007}

\begin{document}
\begin{abstract}
We show that a slight modification of a theorem of Ruzmaikina and Aizenman on competing particle systems on the real line leads to a characterization of Poisson-Dirichlet distributions $PD(\alpha,0)$. 

Precisely, let $\xi$ be a proper random mass-partition i.e. a random sequence $(\xi_i,i\in\N)$ such that $\xi_1\geq \xi_2\geq ...$ and $\sum_i \xi_i=1$ a.s. 
Consider $\{W_i\}_{i\in\N}$, an iid sequence of random positive numbers whose distribution is absolutely continuous with respect to the Lebesgue measure and $\E[W^\lambda]<\infty$ for all $\lambda\in\R$.
It is shown that, if the law of $\xi$ is invariant under the random reshuffling 
$$( \xi_i,i\in\N) \mapsto \left(\frac{\xi_i W_i}{\sum_j \xi_j W_j},i\in\N\right)$$
where the weights are reordered after evolution,
then it must be a mixture of Poisson-Dirichlet distribution $PD(\alpha,0)$, $\alpha\in(0,1)$.
\end{abstract}
\maketitle

\section{The Ruzmaikina-Aizenman Theorem}
Consider $X:=\{X_i\}_{i\in\N}$, a locally finite point process such that $X$ has a maximum almost surely. Plainly, the configuration of points can be ordered decreasingly and we can write $X=(X_i,i\in\N)$ where the indexing corresponds to the decreasing ordering $X_1\geq X_2\geq ...$ Ruzmaikina and Aizenman considered the following evolution of the configuration of points \cite{RA} (see also \cite{liggett} for a similar setup): each point $X_i$ is incremented independently
\begin{equation}
X_i\mapsto X_i +h_i
\label{evolution}
\end{equation}
where the variables $h_i$'s are iid and independent of $X$. Such a dynamics can be seen as a competition within a crowd of points whose evolution is stochastic and uncorrelated. 

In this framework, the law of $X$ is said to be {\it quasi-stationary} if the distribution of the gaps of $X$ is invariant under the map \eqref{evolution}. Namely let $\tilde{X}:=(\tilde{X}_i,i\in\N)=\left(X_i+h_i, i\in\N\right)_\downarrow$ be the evolved point process. Here the symbol $\downarrow$ means that the points are ordered in decreasing order. Quasi-stationarity translates into
$$ \left(\tilde{X}_i-\tilde{X}_{i+1}, i\in\N\right)\distrib\left(X_i-X_{i+1}, i\in\N\right).$$
In a general situation, the configuration of points $\{X_i+h_i\}_{i\in\N}$ might not be locally finite nor have a maximum. If the law of $h$ has a density $g(y)dy$, the laws of point processes for which this property holds are called {\it $g$-regular}. A sufficient condition for $g$-regularity is simply the finiteness of $\E_{X,h}[\#\{i:X_i+h_i\geq y\}]$ for all $y\in\R$ where $\E_{X,h}$ is the expectation over $X$ and $h$. 

Note that quasi-stationary processes must have one point or an infinite number of points almost surely. This is due to the fact that in the case of finite points the gaps between points spread due to the convolution of increments \cite{RA}. The one-point process is trivially quasi-stationary so that we can restrict ourselves to infinite systems. It is easily checked that any Poisson point process with intensity measure $\rho e^{-\rho y}dy$ (and thus any mixture of these) is quasi-stationary for any $h$ with finite exponential moment.
The Ruzmaikina-Aizenman theorem shows that they are the only one up to some assumptions on the density of points. 
\begin{thm1}[Ruzmaikina-Aizenman]
Consider a probability measure with density $g$ on $\R$ such that for all $\lambda\in\R$
\begin{equation}
\int_\R e^{\lambda y}g(y)dy<\infty.
\label{cond g}
\end{equation}
Let $\{h_i\}_{i\in\N}$ be an iid sequence of $g$-distributed random variables and $\mu$ be the law of a point process $X:=(X_i,i\in\N)$ in $\R$ that is almost surely locally finite and has a maximum. Moreover assume that there exist $A>0$ and $r>0$ such that for all $d\geq0$
\begin{equation}
\E_\mu\left[\#\{i\in\N: X_1-X_i\leq d\}\right]\leq Ae^{r d}
\label{tail assump}
\end{equation}

If $\mu$ is quasi-stationary under the evolution \eqref{evolution},
then, as far as the gaps of $X$ are concerned, $\mu$ is a mixture of Poisson point processes with intensity measures $\rho e^{-\rho y}dy$ for $\rho>0$.
\end{thm1}
Interestingly, the quasi-stationarity property holds also when $h$ has a lattice law for Poisson processes with exponential density (see e.g. Proposition 3.1 in \cite{RA}). However, one expects such $h$ not to single out this family. One reason for this is that the Bahadur-Rao theorem, used to get the exponential form of the density in the proof of Theorem 1, takes a different form for lattice laws (see e.g. \cite{dembo} Theorem 3.7.4).

For the proof of Theorem 1, the authors chose without loss of generality that $X_1=0$ for all realizations of $X$. One consequence of this choice is that the condition \eqref{tail assump} ensures the $g$-regularity of $\mu$ when $\mu$ is quasi-stationary. We will see that a different choice of normalization of the initial configuration $X$ leads to a modified version of the theorem.

\subsection{A variation of the theorem}
Let us look at a different class of point processes. Consider $X=\{X_i\}_{i\in\N}$ for which there exists a $\beta>0$ such that
\begin{equation}
\sum_{i\in \N} e^{\beta X_i}<\infty \text{ a.s.}
\label{cond sum}
\end{equation}
Note that the configuration of points of such processes are automatically locally finite and have a maximum. As before, they can be ordered in decreasing order and we write $X=(X_i,i\in\N)$. Point processes of this kind are natural in a physics setting. If one considers $X$ as a collection of random energies, then this class corresponds to random systems with almost sure finite partition function  at inverse temperature $\beta$.

Again we consider the evolution \eqref{evolution} where $h$ has law $g(y)dy$.
One advantage of the point processes we just introduced is the fact that their law is always $g$-regular for $g$ satisfying \eqref{cond g}. Indeed, to prove that the evolved configuration is locally finite and has a maximum, it suffices to show that
$$ \sum_{i\in\N} e^{\beta(X_i+h_i)}<\infty \text{ a.s.}$$
The above would follow if $\E_{\mu,h}\left[\frac{\sum_{i\in\N} e^{\beta(X_i+h_i)}}{\sum_{i\in\N} e^{\beta X_i}}\right]<\infty$ a.s. But this is clear by the fact that $h$ is independent of $X$ and $g$ has finite Laplace transform at $\beta$. 

Is the condition \eqref{tail assump} necessary for the conclusion of Theorem 1 to hold? Roughly speaking, such a condition yields a control on the tail of the density of points. However, condition \eqref{cond sum} is also in some sense a bound on the density of points. It turns out that a slight modification of an estimate needed in the proof of Theorem 1 shows that condition \eqref{cond sum} can effectively replace the assumption \eqref{tail assump}. The details of the argument are given in the next section. 
\begin{thm1}[Modified version]
Let $g$ be as in Theorem 1.
Let $\mu$ be the of the law of a point process $X$ satisfying equation \eqref{cond sum}.

If $\mu$ is quasi-stationary under the evolution \eqref{evolution},
then, as far as the gaps of $X$ are concerned, $\mu$ is a mixture of Poisson point processes with intensity measures $\rho e^{-\rho y}dy$ for $\rho>0$.
\label{ra mod}
\end{thm1}

An important feature of processes satisfying \eqref{cond sum} is the fact that they are closely related to random mass-partitions. A {\it random mass-partition} is a random sequence $(\xi_i,i\in\N)$ such that $\xi_1\geq \xi_2\geq ...$ and $\sum_i \xi_i\leq 1$ It is said to be {\it proper} if $\sum_i \xi_i=1$ a.s. We will say that a random mass-partition is {\it trivial} if $\xi_1=1$ a.s. If $X$ satisfies \eqref{cond sum}, then we can define the $X$-measurable random mass-partition
\begin{equation}
\xi:=\left(\frac{e^{\beta X_i}}{\sum_j e^{\beta X_j}},i\in\N\right).
\label{rmp X}
\end{equation}
Conversely, if $\xi$ is a random mass-partition, then the point process $X:=\left(\log \xi_i,i\in\N\right)$ satisfies \eqref{cond sum} for $\beta=1$. 

Important examples of random mass-partitions are provided by the so-called Poisson-Kingman partitions. They are constructed in the following way. Let $\Lambda$ be a measure on $\R^+$ such that $\int_0^\infty(1\wedge y)\Lambda(dy)<\infty$ and $\Lambda(\R^+)=\infty$. Consider $\eta$, a Poisson point process with intensity measure $\Lambda$. The above conditions are necessary and sufficient for $\eta$ to have an infinite number of atoms and for which $\sum_i \eta_i<\infty$ almost surely. A Poisson-Kingman partition is obtained by simply taking $\xi_i:=\eta_i/\sum_j \eta_j$. The Poisson-Dirichlet distributions $PD(\alpha,0)$ is a subclass of Poisson-Kingman partitions obtained from measures $\Lambda$ of the form $\alpha s^{-\alpha-1} ds$, $\alpha\in (0,1)$. Note that these are exactly the image measures by an exponential map that sends $\R$ to $\R^+$ of the measures with exponential density. Poisson-Dirichlet distributions are fundamental objects in several areas of probability such as coagulation, fragmentation and random partitions (see \cite{bertoin} for details).

In the framework of random mass-partition, the stochastic map \eqref{evolution} simply becomes a multiplicative shift followed by a renormalization
\begin{equation}
\xi_i\mapsto \frac{\xi_iW_i}{\sum_j\xi_j W_j}
\label{evolution rmp}
\end{equation}
where $W_i=e^{\beta h_i}$.
The weights of $\xi$ are reshuffled under the above map as were the points of the process $X$ and are reordered decreasingly after evolution.
It is easy to see that the law of $X$ is quasi-stationary if and only if the law of $\xi$ constructed from $X$ in equation \eqref{rmp X} is invariant under the evolution \eqref{evolution rmp} i.e.
\begin{equation*}
\left(\frac{\xi_iW_i}{\sum_j\xi_jW_j},i\in\N\right)_\downarrow\distrib \xi
\label{invar rmp}
\end{equation*}
Note that the trivial random mass-partition is plainly invariant under the stochastic map \eqref{evolution rmp}. It can be easily checked using the quasi-stationarity of Poisson process with intensity measures $\rho e^{-\rho y}dy$ that Poisson-Dirichlet laws $PD(\alpha,0)$ (and any mixture of these) are non-trivial examples of stationary laws under \eqref{evolution rmp}. Our goal is to show that they are the only ones.

We would like to answer this question by using the Ruzmaikina-Aizenman Theorem. If we want to apply the original version to the point process associated to $\xi$, Assumption \eqref{tail assump} must be satisfied. In terms of mass-partition, it translates into
$\E_\xi\left[\#\{i: \xi_i/\xi_1 \geq \delta \}\right]\leq A \delta^{-\lambda}
$
for all $\delta\in(0,1)$ and for some $\lambda>0$ and $A>0$.
A sufficient condition (using Markov's inequality) for the above to hold is the finiteness for some $\lambda>0$ of $
\E_\xi\left[\sum_{i\in\N}\left(\xi_i/\xi_1\right)^\lambda\right].$
One can construct examples of random mass-partitions that do not satisfy the above.
However, the modified version of Theorem 1 applies to any proper random mass partition thereby yielding the desired characterization of the Poisson-Dirichlet distributions $PD(\alpha,0)$, $\alpha\in (0,1)$.
\begin{thm2}
Let $\xi$ be a non-trivial proper random mass-partition. Let $\{W_i\}_{i\in\N}$ be an iid sequence of positive numbers whose distribution has a density on $\R^+$ and $\E[W^\lambda]<\infty$ for any $\lambda\in\R$.

If $\xi$ is invariant under the evolution \eqref{evolution rmp}, then its law is a mixture of Poisson-Dirichlet distributions $PD(\alpha,0)$, $\alpha\in (0,1)$.
\end{thm2}
\begin{proof}
The law of $\xi$ induces a point process on $\R$, namely $X:=\left(\log \xi_i,i\in\N\right)$, which plainly respects the summability condition \eqref{cond sum}. Moreover, $\log W$ satisfies the conditions on $h$ in Theorem 1.
Therefore we can use the modified version of Theorem 1 to conclude that $\{X_i\}_{i\in\N}$ is a mixture of Poisson processes with intensity measure $\alpha e^{-\alpha s}ds$ for $\alpha>0$ as far as the gaps are concerned. 
Without loss of generality, we can assume that $X$ is actually a mixture of Poisson processes because the law of the random mass-partition depends only on the distribution of the gaps. 
One can verify that the process $(e^{X_i},i\in\N)$ must therefore be a mixture of Poisson process on $\R^+$ with intensity measure $\alpha s^{-\alpha-1}ds$. However, the positions of the points of this process are summable if and only if $\alpha\in(0,1)$. By definition of Poisson-Dirichlet distributions, the normalized process $\xi$ is a mixture of random mass partitions $PD(\alpha,0)$.
\end{proof}
It is important to stress that every law $PD(x,0)$ is invariant under the evolution for any $W$ satisfying the assumptions. Therefore, nothing can be said on the specific composition of the mixture constituting $\xi$. However, if $\xi$ is ergodic under the evolution, its law must be pure.

It would be of interest to find an alternative proof of the above result without appealing to the point process associated to $\xi$. Also one might wonder if the extended two-parameter family of Poisson-Dirichlet variables $PD(\alpha,\theta)$ is characterized by an invariance property under a general stochastic transformation in the flavor of Theorem 2. See \cite{bertoin} for the definition of $PD(\alpha,\theta)$ variables and \cite{diaconis} for a characterization result in the case $PD(0,1)$.

\section{Modification of the proof}
The original proof of the Ruzmaikina-Aizenman Theorem only needs to be slightly modified for Theorem 1 to hold in its modified version. Roughly speaking, one needs to retrieve bounds on the displacement of the front of the crowd of points under the new assumption. In this section, we first give a short summary of the original proof of the result. Secondly, we prove the necessary estimates for the result to hold under the modified assumption of Theorem 1.
\subsection{A summary of the original proof}
Quasi-stationarity is by definition a property that only depends on the law of the gaps of the point process $X$. Therefore, to characterize the quasi-stationary laws it suffices to find an explicit expression for the generating functional
\begin{equation}
G_\mu(f)=\E_\mu\left[\exp \sum_i f(X_1-X_i)\right]
\label{G mu}
\end{equation}
where $f$ is a measurable positive compactly supported function.

Let $\mu$ be a Poisson point process with intensity measure $-dF$ where $F$ is a smooth monotone decreasing function with $F(y)\to 0$ when $y\to\infty$. 
In the case where $\mu$ is a Poisson point process with intensity measure $-dF$, it is not hard to check that this generating functional has the form
\begin{equation}
G_{PP(-dF)}(f)=\int_\R d(e^{-F(x)}) \exp\left(-\int_{(\infty,x]}(1-e^{-f(y)})-dF(y)\right).
\label{G PP}
\end{equation}

To prove Theorem 1, the authors show that if $\mu$ is quasi-stationary then its generating functional \eqref{G mu} is given by an integral over a parameter $\rho>0$ of the functional \eqref{G PP} with $F(y)=e^{-\rho y}$. 

The crucial point in the proposed modification of the original proof is to notice that the generating functional \eqref{G mu} is by definition independent of a global shift of the initial configuration of points $X$ (see \cite{RA} p. 88). Therefore, one has the freedom to choose such a shift. The convenient choice made in \cite{RA} is a shift by the leader i.e. one assumes that almost surely $X_1=0$ (corresponding to the shift $X_i\mapsto X_i- X_1$). The modified version of the theorem is based on the shift by the tail namely $X_i\mapsto X_i-\frac{1}{\beta} \log \sum_i e^{\beta X_i}$ so that one can assume $\sum_i e^{\beta X_i}=1$ a.s.

No matter what global shift of the initial configuration is chosen, the first part of the proof holds (Section 5 in \cite{RA}). Precisely, one defines the function
$$ F_{X,\tau}(y):=\sum_i \prob_h(S_i(\tau)+X_i\geq y)$$
for $S_i(\tau):=\sum_{t=1}^\tau h_i(t)$ where the $h_i(t)$'s are iid copies of $h_i$ and $\prob_h$ is the integration over the increments given $X$. The function $F_{X,\tau}$ represents the expected number of points beyond $y$ after $\tau$ steps conditionally on the initial configuration $X$. One also defines $Z_{X,\tau}$ as the point where
$$ F_{X,\tau}(Z_{X,\tau})=1.$$
The random variable $Z_{X,\tau}$ can be interpreted as the expected position of the leading edge of the crowd of points giving $X$. A normalizing shift $\mathcal{N}$ is also introduced 
$$\mathcal{N}F_{X,\tau}(y):=F_{X,\tau}(y+Z_{X,\tau})$$ 
so that $\mathcal{N}F_{X,\tau}(0)=1$ for all $X$. Clearly the generating functional of the gaps is invariant under the normalizing operation. The first result of \cite{RA} is to show that a quasi-stationary measure must be a mixture of Poisson process with intensity measures $-d\mathcal{N}F_{X,\tau}(y)$.
\begin{thm5.1}[\cite{RA}]
Let $\mu$ be a $g$-regular quasi-stationary measure under the evolution \eqref{evolution}. Then for any measurable positive function $f$ with compact support
$$G_\mu(f)=\lim_{\tau\to\infty}\E_\mu\left[G_{PP(-d\mathcal{N}F_{X,\tau})} (f)\right]. $$
\end{thm5.1}
The main idea of the proof is to show that a quasi-stationary point process tends to a mixture of Poisson process due to the uncorrelated nature of the evolution.
Recall that in the case where $X$ satisfies the summability condition \eqref{cond sum}, its law is automatically $g$-regular so that we can drop this assumption in the above theorem.

The second part of the proof (Section 6 in \cite{RA}) is where the assumption \eqref{tail assump} is needed. One would like to prove that the Poisson processes have an intensity measure given by the Laplace transform of a finite measure.
\begin{thm6.1}[\cite{RA}]
Under the assumption of Theorem 1, there exists a probability measure $\nu(d\rho)$ on the space $\mathcal{M}$ of finite measures on $\R$ such that for any measurable positive compactly supported function $f$
$$G_\mu(f)=\int_{\mathcal{M}}\nu(d\rho)\text{ }G_{PP(-dR_\rho)}(f) $$
where $R_\rho(y)=\int_0^\infty e^{-yu}\rho(du)$. 
\end{thm6.1}
The proof of the above is based on a lemma that gives control on the random variable $Z_{X,\tau}$ as well as on the number of points close to the leader at time $\tau$. The assumption \eqref{tail assump} of Theorem 1 as well as the choice $X_1=0$ are needed in the proof of this lemma only. Therefore, if the same control can be achieved through a different choice of initial shift and assumption on $X$, the result would still hold. We state the lemma in its original form. Its proof under the assumption of the modified version of Theorem 1 is given in the next subsection.
\begin{lem6.2}[\cite{RA}]
Let $\mu$ be a $g$-regular quasi-stationary measure that satisfies the assumption of Theorem 1. Assume also that initially $X_1=0$ almost surely. Then the following hold
\begin{enumerate}
\item For any $\epsilon>0$, there exists $\tau$ large enough and $\lambda>0$ such that on a subset of $\mu$-measure $1-\epsilon$
$$ Z_{X,\tau}\leq \frac{S}{2\lambda}\tau +\text{const}$$
where $S:=\log \int_\R e^{2\lambda y}g(y)dy$.

\item Let $A_{\tau,d,K,M}$ be the event that the configurations obtained after $\tau$ steps will not have more than $M$ points within distance $d$ of the leading point and all of them made a total jump less than $K\tau +Z_{X,\tau}-X_i$ from time $0$ to $\tau$. Then 
$$ \prob_{\mu,h}(A^c_{\tau,d,K,M})\leq \Gamma_1(M,d)+\Gamma_2(\tau) + C e^{-\delta(K-K_0)\tau}$$
where $\delta,K_0,C>0$ and $\Gamma_1(M,d)\to 0$ as $M\to\infty$ for fixed $d$ and $\Gamma_2(\tau)\to 0$ as $\tau\to\infty$.
\end{enumerate}
\end{lem6.2}

Finally the last part of the proof (Section 7 in \cite{RA}) consists in showing that the functions $R_\rho$ are almost surely pure exponentials. The argument is based on general monotonicity properties of $R_\rho$ under convolution that holds in a general setting.
\subsection{The modified estimate}
As argued in the last section, the modified version of Theorem 1 will be established if Lemma 6.2 is proven under the new hypothesis. Precisely, we now assume that the initial configuration is reshifted by $\frac{1}{\beta}\log \sum_i e^{\beta X_i}$ so that $\sum_ie^{\beta X_i}=1$ a.s. 
\begin{lem}
If $X$ is a point process such that there exists $\beta>0$ for which $\sum_i e^{\beta X_i}=1$ a.s., then the conclusion of Lemma 6.2 holds.
\end{lem}
\begin{proof}[Proof of (1)]
We show that 
$$ F_{X,\tau}(y)\leq e^{v_\beta\tau-\beta y} \text{ a.s.}$$
for some $v_\beta>0$. 
This will mean that $Z_{X,\tau}\leq \frac{v_\beta}{\beta}\tau$ almost surely.
The above is clear from a simple application of Markov's inequality with the function $e^{\beta y}$. We have
\begin{align*}
F_{X,\tau}(y)&=\sum_i \prob_h(S_i(\tau) +X_i\geq y)\leq e^{-\beta y}\sum_i e^{\beta X_i}\E_h\left[e^{\beta S_i(\tau)}\right]= e^{v_\beta\tau-\beta y}
\end{align*}
where $v_\beta:=\log \E[e^{\beta h}]$ and we have used the fact that $\sum_i e^{\beta X_i}=1$. Note that we actually achieved a better control on $Z_{X,\tau}$ than the one in the original proof as the bound holds almost surely and not only in probability.
\end{proof}
\begin{proof}[Proof of (2)]
Let $C<\E[h]$. Consider $B_{\tau,K,C}$ the event that at least one point made a jump in $\tau$ steps greater than $-X_i+(C+K)\tau$. Then the probability of the complement of the event $A_{\tau,d,K,M}$ is bounded above by
$$
\prob_{\mu,h}(A^c_{\tau,d,K,M})\leq \prob_{\mu,h}\left(\#\{i:\tilde{X}_1(\tau)-\tilde{X}_i(\tau)\leq d\}\geq M\right)+\prob_{\mu,h}(B_{\tau,K,C})+ \prob_\mu(Z_{X,\tau}\leq C\tau) 
$$
where $\tilde{X}(\tau)$ is the ordered process at time $\tau$.
By quasi-stationarity, the first term equals the probability at time $0$ that the number of points within distance $d$ of the leader is greater than M. But as the configuration is locally finite, it goes to $0$ as $M\to\infty$ for fixed $d$.
As for the second term, we have by a similar argument as the above proof
\begin{align*}
 \prob_{\mu,h}(B_{\tau,K,C})&=\prob_{\mu,h}\left(\bigcup_i\left\{S_i(\tau)+X_i\geq (C+K)\tau\right\}\right)\\
 &\leq \sum_i \prob_{\mu,h}\left(\left\{S_i(\tau)+X_i\geq (C+K)\tau\right\}\right)\\
 &\leq \sum_i e^{-\beta (C+K)\tau}\E_{\mu,h}\left[e^{\beta X_i+\beta S_i(\tau)}\right]\\
 &= e^{-\tau((C+K)\beta -v_\beta)}
\end{align*}
where we have used the fact that $X$ is independent of $h$.
As $\beta$ is fixed, we can choose $K$ large enough so that the desired asymptotics holds. 

It remains to estimate the third term. Choose $C:=\E[h]-1$. Let $Z'_{X,\tau}$ be the position of the front for the original choice of initial shift $X_1=0$. By definition, $Z_{X,\tau}=Z'_{X,\tau}-r_X$ where $r_X:=\frac{1}{\beta}\log \sum_i e^{\beta X_i}- X_1$. Note that $r_X>0$. The asymptotics of $Z'_{X,\tau}$ is proven in \cite{RA} (see p.99 following equation 6.11)  for all $C'<\E[h]$ and some function $\Gamma'_{C'}(\tau)$.
\begin{equation}
\prob_\mu(Z'_{X,\tau}\leq C'\tau)\leq \Gamma'_{C'}(\tau)\to 0 \text{ as $\tau\to\infty$}.
\label{eqn z}
\end{equation}
Equation \eqref{eqn z} induces the asymptotics of $Z'_{X,\tau}$
\begin{align*}
 \prob_\mu(Z_{X,\tau}\leq C\tau)&= \prob_\mu(Z'_{X,\tau}\leq C\tau+r_X)\\
 &=\prob_\mu(Z'_{X,\tau}\leq \tau(C+r_X/\tau), \text{ } r_X/\tau>\delta)+\prob_\mu(Z'_{X,\tau}\leq \tau(C+r_X/\tau), \text{ } r_X/\tau\leq \delta)\\
& \leq \prob_\mu(r_X/\tau>\delta)+\prob_\mu(Z'_{X,\tau}\leq \tau(C+\delta)).
\end{align*}
We conclude that $\prob_\mu(Z_{X,\tau}\leq C\tau)\to 0$ as $\tau\to \infty$ if we choose any $0<\delta <1$ using equation \eqref{eqn z} and the fact that $r_X$ is independent of $\tau$.
\end{proof}



%

\end{document}